\documentclass{article}
\usepackage{latexsym}
\usepackage{amssymb}
\usepackage{amsmath}
\usepackage{mathrsfs}

\allowdisplaybreaks

\newtheorem{theorem}{Theorem}[section]
\newtheorem{lemma}{Lemma}[section]
\newtheorem{proposition}{Proposition}[section]
\newtheorem{definition}{Definition}[section]
\newtheorem{remark}{Remark}[section]

\newtheorem{example}{Example}[section]

\begin{document}

\title{\bf Attractor-repeller pair, Morse decomposition and Lyapunov
function for random dynamical systems\footnote{This work is
supported by the 985 project of Jilin University and Graduate
Innovation Lab of Jilin University.}}

\author{Zhenxin Liu\footnote{Corresponding author.}, Shuguan Ji and Menglong Su\\
{\small College of Mathematics, Jilin University,  Changchun 130012,
People's Republic of China}\\
{\small zxliu@email.jlu.edu.cn; jesenliu@163.com}}

\maketitle

\begin{abstract}
In the stability theory of dynamical systems, Lyapunov functions
play a fundamental role. In this paper, we study the
attractor-repeller pair decomposition and Morse decomposition for
compact metric space in the random setting. In contrast to
\cite{Cra}, by introducing slightly stronger definitions of random
attractor and repeller, we characterize attractor-repeller pair
decompositions and Morse decompositions for random dynamical systems
through the existence of Lyapunov functions. These
characterizations, we think, deserve to be known widely.
\end{abstract}

\section{Introduction and main result}

In the stability theory of dynamical systems, Lyapunov functions
have been playing a fundamental role ever since first introduced in
Lyapunov's 1892 thesis \cite{Lia}. The simple idea of linking
dynamics and topology by means of functions decreasing along
trajectories has subsequently been instrumental in the development
of sophisticated tools such as e.g. Morse decompositions and Floer
homology. Attractor-repeller pairs and Morse sets are special
invariant sets that plays an important role in understanding the
asymptotic behavior of a topological dynamical system defined on a
compact metric state space. A complete treatment of Morse theory for
deterministic case can be found in the monograph of Conley
\cite{Con}. Among interesting results in \cite{Con} is a proposition
which claims that mutually disjoint invariant sets are an
attractor-repeller pair if and only if there exists a Lyapunov
function for those sets. The result is then developed in Huang
\cite{Hua} as a criterion to detect Morse sets and Morse
decomposition once a Lyapunov function for those sets can be
constructed.

Random dynamical systems (RDS) arise in the modeling of many
phenomena in physics, biology, economics, climatology, etc and the
random effects often reflect intrinsic properties of these phenomena
rather than just to compensate for the defects in deterministic
models. The history of study of random dynamical systems goes back
to Ulam and von Neumann \cite{UV} and it has flourished since the
1980s due to the discovery that the solutions of stochastic ordinary
differential equations yield a cocycle over a metric dynamical
system which models randomness, i.e. a random dynamical system. In
developing a comprehensive theory of random dynamical systems,
members of the Bremen Group started establishing analogous notions,
techniques and results for the stochastic setting. Lyapunov
functions for RDS were introduced by Arnold and Schmalfuss
\cite{Ar2}, and Crauel et al \cite{Cra} established Morse
decompositions and studied some of their basic properties for RDS.

The present paper contributes to this ongoing process. We study the
attractor-repeller pair decomposition and Morse decomposition for
compact metric space in the random setting. In contrast to
\cite{Cra}, by introducing slightly stronger definitions of random
attractor and repeller, we can construct {\em measurable} Lyapunov
functions for attractor-repeller pair and Morse decomposition for
RDS. And moreover we also prove that the existence of {\em
continuous} Lyapunov functions is also the sufficient condition to
conclude that two (or finite) mutually disjoint invariant random
compact sets constitute an attractor-repeller pair (or a Morse
decomposition) for RDS.

The paper is organized as follows. In Section 2, we recall some
basic definitions and results for RDS. In Section 3, we give the
definitions of limit set (omega-limit set and alpha-limit set),
attractor and repeller for RDS. In Section 4, we study the
attractor-repeller pair decomposition on compact metric space in the
random setting and characterize it by Lyapunov function. And at last
we characterize the Morse decompositions on compact metric space
through the Lyapunov function and give a simple example in Section
5.

\section{Random dynamical systems}

Throughout the paper all assertions about $\omega$ are assumed to
hold on a $\theta$ invariant set of full measure unless otherwise
stated. First we give the definition of continuous random dynamical
systems (cf. Arnold \cite{Ar1}).

\begin{definition}\rm
Let $X$ be a metric space with a metric $d$. A (continuous) {\em
random dynamical system (RDS)}, shortly denoted by $\varphi$,
consists of two ingredients:
\begin{itemize}
    \item A model of
the noise, namely a metric dynamical system $(\Omega, \mathscr F,
\mathbb P, (\theta_t)_{t\in \mathbb R})$, where $(\Omega, \mathscr
F, \mathbb P)$ is a probability space and $(t,\omega)\mapsto
\theta_t\omega$ is a measurable flow which leaves $\mathbb P$
invariant, i.e. $\theta_t\mathbb P=\mathbb P$ for all $t\in \mathbb
R$. For simplicity we also assume that $\theta$ is ergodic under
$\mathbb P$, meaning that a $\theta$-invariant set has probability 0
or 1.
    \item  A model of the system perturbed by noise,
namely a cocycle $\varphi$ over $\theta$, i.e. a measurable mapping
$\varphi: \mathbb R\times \Omega\times X \rightarrow X,
(t,\omega,x)\mapsto\varphi(t,\omega,x)$, such that
$(t,x)\mapsto\varphi(t,\omega,x)$ is continuous for all
$\omega\in\Omega$ and the family
$\varphi(t,\omega,\cdot)=\varphi(t,\omega):X\rightarrow X$ of random
self-mappings of $X$ satisfies the cocycle property:
\begin{equation}\label{phi}
\varphi(0,\omega)={\rm id}_X, \varphi(t+s,\omega)=\varphi(t,\theta_s
\omega)\circ\varphi(s,\omega)\quad {\rm for ~all}\quad t,s\in\mathbb
R,\omega\in\Omega.
\end{equation}
\end{itemize}
\end{definition}

It follows from (\ref{phi}) that $\varphi(t,\omega)$ is a
homeomorphism of $X$, and the fact
\[\varphi(t,\omega)^{-1}=\varphi(-t,\theta_t\omega)\]
is very useful in the following.

Any mapping from $\Omega$ into the collection of all subsets of $X$
is said to be a {\it multifunction} (or a set valued mapping) from
$\Omega$ into X. We now give the definition of random set, which is
a fundamental concept for RDS.

\begin{definition}\rm
Let $X$ be a metric space with a metric $d$. The multifunction
$\omega\mapsto D(\omega)$ taking values in the closed/compact
subsets of $X$ is said to be a {\em random closed/compact set} if
the mapping $\omega\mapsto {\rm dist}_X(x,D(\omega))$ is measurable
for any $x \in X$, where ${\rm dist}_X(x,B):=\inf_{y\in B}d(x,y)$.
The multifunction $\omega\mapsto U(\omega)$ taking values in the
open subsets of $X$ is said to be a {\em random open set} if
$\omega\mapsto U^c(\omega)$ is a random closed set, where $U^c$
denotes the complement of $U$.
\end{definition}

Afterwards, we also call a multifunction $D(\omega)$ measurable for
convenience if the mapping $\omega\mapsto {\rm dist}_X(x,D(\omega))$
is measurable for any $x \in X$.

\begin{definition}\rm
A random set $D(\omega)$ is said to be {\em forward invariant} under
the RDS $\varphi$ if $\varphi(t,\omega)D(\omega)\subset
D(\theta_t\omega)$ for all $t\ge 0$; It is said to be {\em backward
invariant} if $\varphi(t,\omega)D(\omega)\supset D(\theta_t\omega)$
for all $t\ge 0$; It is said to be {\em invariant} if
$\varphi(t,\omega)D(\omega)=D(\theta_t\omega)$ for all $t\in\mathbb
R$.
\end{definition}

Now we enumerate some basic results about random sets in the
following proposition, for details the reader can refer to Castaing
and Valadier \cite{Cas}, Crauel \cite{Cr} and Arnold \cite{Ar1} for
instance.
\begin{proposition}\label{set}
Let X be a Polish space, then the following assertions hold: \\
(i) if D is a random closed set, then so is the  closure of $D^c$;\\
(ii) if D is a random open set, then the closure $\overline{D}$ of
$D$
is a random closed set;\\
(iii) if D is a random closed set, then {\rm int}$D$, the interior
of $D$, is a random open set; \\
(iv) if $\{D_n, n \in\mathbb N\}$ is a sequence of random closed
sets and there exists $n_0 \in\mathbb N$ such that $D_{n_0}$ is a
random compact set, then
$\bigcap_{n\in\mathbb N}D_n$ is a random compact set;\\
(v) if $f:\Omega\times X\rightarrow X$ is a function such that
$f(\omega,\cdot)$ is continuous for all $\omega$ and $f(\cdot,x)$ is
measurable for all $x$, then $\omega\mapsto f(\omega,D(\omega))$ is
a random compact set provided that $D(\omega)$ is a random compact
set.
\end{proposition}

\section{Limit set, attractor and repeller}

\begin{definition}\rm
For any given random set $D(\omega)$, we denote $\Omega_D(\omega)$
the {\em omega-limit set of $D(\omega)$}, which is determined as
follows:
\[
\Omega_D(\omega):=\bigcap_{T\ge 0}\overline{\bigcup_{t\ge
T}\phi(t,\theta_{-t}\omega)D(\theta_{-t}\omega)};
\]
and we denote $\alpha_D(\omega)$ the {\em alpha-limit set of
$D(\omega)$}, which is determined as follows:
\[
\alpha_D(\omega):=\bigcap_{T\ge 0}\overline{\bigcup_{t\ge
T}\phi(-t,\theta_{t}\omega)D(\theta_{t}\omega)}.
\]
\end{definition}

\begin{definition}\rm
For given two random sets $D(\omega),A(\omega)$, we say $A(\omega)$
(pull-back) {\em attracts} ({\em repels}) $D(\omega)$ if
\[
\lim_{t\rightarrow\infty}d(\varphi(t,\theta_{-t}\omega)D(\theta_{-t}\omega)|A(\omega))=0
~(\lim_{t\rightarrow-\infty}d(\varphi(t,\theta_{-t}\omega)D(\theta_{-t}\omega)|A(\omega))=0)
\]
 holds almost surely, where $d(A|B)$ stands for the Hausdorff
semi-metric between two sets $A,B$, i.e. $d(A|B):={\rm sup}_{x\in
A}{\rm inf}_{y\in B}d(x,y)$.
\end{definition}

\begin{remark}\label{re}\rm
It is well-known that $x\in \Omega_{D}(\omega)$ if and only if
$\exists t_n\rightarrow\infty$, $x_n\in D(\theta_{-t_n}\omega)$ such
that $\varphi(t_n,\theta_{-t_n}\omega)x_n\rightarrow x$,
$n\rightarrow\infty$. If a non-void  random set $D(\omega)$ is
attracted by a random compact set $K(\omega)$, then
$\Omega_D(\omega)\neq\emptyset$ almost surely and it is invariant.
Moreover, $\Omega_D(\omega)$ pull-back attracts $D(\omega)$.
\end{remark}

\begin{definition}\label{def}\rm
(i) An invariant random compact set $A(\omega)$ is called an {\em
(local) attractor} if there exists a random closed neighborhood
$N(\omega)$ of $A(\omega)$ such that $A(\omega)=\Omega_N(\omega)$.
The closed neighborhood $N(\omega)$ is called a {\em fundamental
neighborhood} of $A(\omega)$.\\
(ii) An invariant random compact set $R(\omega)$ is called a {\em
(local) repeller} if there exists a random closed neighborhood
$N(\omega)$ of $R(\omega)$ such that $R(\omega)=\alpha_D(\omega)$.
The closed neighborhood $N(\omega)$ is called a {\em fundamental
neighborhood} of $R(\omega)$.
\end{definition}

The following definition of basin is adopted in \cite{Cra,Liu}.

\begin{definition}\rm
(i) Assume $A(\omega)$ is an attractor with a fundamental
neighborhood $N(\omega)$. Then we call
\[
B(A)(\omega):=\{x|~\varphi(t,\omega)x\in {\rm
int}N(\theta_t\omega)~{\rm for~ some}~t\ge 0\}
\]
the {\em basin of attraction} of $A(\omega)$;\\
(ii) Assume $R(\omega)$ is a repeller with a fundamental
neighborhood $N(\omega)$. Then we call
\[
B(R)(\omega):=\{x|~\varphi(t,\omega)x\in {\rm
int}N(\theta_t\omega)~{\rm for~ some}~t\le 0\}
\]
the {\em basin of repulsion} of $R(\omega)$.
\end{definition}

\begin{remark}\rm
The basins of attractor and repeller are well defined, i.e. they do
not depend on the choice of their fundamental neighborhoods. The
readers can refer to \cite{Cra,Liu} for details.
\end{remark}

\begin{lemma}\label{lem}
Assume $N(\omega)$ is a random closed set and an invariant random
compact set $A(\omega)\subset {\rm int}N(\omega)$ satisfying that
$\Omega_N(\omega)=A(\omega)$, then there exists a forward invariant
random closed set $\hat N(\omega)$ with the same properties as
$N(\omega)$.
\end{lemma}
\noindent{\bf Proof.} Let
\[
\tilde N(\omega):=\overline{\bigcup_{t\ge
0}\varphi(t,\theta_{-t}\omega)N(\theta_{-t}\omega)},
\]
then by Proposition 1.5.1 of \cite{Chu} we have $\tilde N(\omega)$
is a universally measurable forward invariant random closed set and
$A(\omega)\subset{\rm int}\tilde N(\omega)$ (note that
$N(\omega)\subset\tilde N(\omega)$). Now we show that
$\Omega_{\tilde N}(\omega)=A(\omega)$.
\begin{align}
\Omega_{\tilde N}(\omega)&=\bigcap_{T\ge 0}\overline{\bigcup_{s\ge
T}\varphi(s,\theta_{-s}\omega)\tilde N(\theta_{-s}\omega)}\nonumber\\
&=\bigcap_{T\ge 0}\overline{\bigcup_{s\ge
T}[\varphi(s,\theta_{-s}\omega)\bigcup_{t\ge
0}\varphi(t,\theta_{-t}\circ\theta_{-s}\omega)N(\theta_{-t}\circ\theta_{-s}\omega)]}\nonumber\\
&=\bigcap_{T\ge 0}\overline{\bigcup_{s\ge T}\bigcup_{t\ge
0}\varphi(s,\theta_{-s}\omega)\circ\varphi(t,\theta_{-t}\circ\theta_{-s}\omega)N(\theta_{-t}\circ\theta_{-s}\omega)}
\nonumber\\
&=\bigcap_{T\ge 0}\overline{\bigcup_{s\ge T}\bigcup_{t\ge
0}\varphi(s+t,\theta_{-s-t}\omega)N(\theta_{-s-t}\omega)}
\nonumber\\
&=\bigcap_{T\ge 0}\overline{\bigcup_{s\ge
T}\varphi(s,\theta_{-s}\omega)N(\theta_{-s}\omega)} \nonumber\\
&=\Omega_N(\omega)=A(\omega), \nonumber
\end{align}
where the second ``=" holds since for any random set $D(\omega)$ we
have
\[ \overline{\bigcup_{t\ge
T}\varphi(t,\theta_{-t}\omega)D(\theta_{-t}\omega)}=
\overline{\bigcup_{t\ge
T}\varphi(t,\theta_{-t}\omega)\overline{D(\theta_{-t}\omega)}}.
\]
By Lemma 2.7 in \cite{Cr}, there exists an $\mathscr F$-measurable
random closed set $\hat N(\omega)=\tilde N(\omega)$ almost surely.
This completes the proof of the lemma. \hfill$\Box$

\begin{remark}\rm
By Lemmas \ref{lem}, for a given attractor, we can always choose a
forward invariant random closed set as its fundamental neighborhood.
Hence from now on when we talk about fundamental neighborhood we
mean a forward invariant one; when we say ``{\em strong}"
fundamental neighborhood $N(\omega)$ of $A(\omega)$ we mean that
$N(\omega)$ is a forward invariant fundamental neighborhood and it
satisfies that $\varphi(t,\omega)x\in{\rm int}N(\theta_t\omega)$ for
arbitrary $x\in N(\omega)$ and  $t>0$.
\end{remark}

\section{Attractor-repeller pair and Lyapunov function}

From now on we assume $X$ is a compact metric space, i.e. we will
study the attractor-repeller pair decomposition and Morse
decomposition on compact metric space.

In this section, we mainly consider the relation between
attractor-repeller pair and Lyapunov function for RDS.

\begin{lemma}
Assume $\varphi$ is an RDS on a compact metric space $X$,
$A(\omega)$ is an attractor of $\varphi$ with a fundamental
neighborhood $N(\omega)$ and the basin of attraction $B(A)(\omega)$.
Then $R(\omega):=X\backslash B(A)(\omega)$ is a random repeller with
the basin of repulsion $X\backslash A(\omega)$ and a fundamental
neighborhood $X\backslash{\rm int}N(\omega)$.
\end{lemma}
\noindent{\bf Proof.} Since $N(\omega)$ is a forward invariant
random compact set, we have that $\hat N(\omega):=X\backslash{\rm
int}N(\omega)$ is a backward invariant random compact set (see page
35 of \cite{Ar1}). Denote $\hat R(\omega):=\alpha_{\hat N}(\omega)$.
Then $\hat R(\omega)$ is a random repeller with a fundamental
neighborhood $\hat N(\omega)$. By the definition of alpha-limit set,
the facts $R(\omega)\subset\hat N(\omega)$ and the invariance of
$R(\omega)$ we have $R(\omega)\subset\hat R(\omega)$. If there
exists $x_0\in\hat R(\omega)\backslash R(\omega)$, then $x_0\in
B(A)(\omega)$. Therefore there exists some $t_0\ge 0$ such that
$\varphi(t_0,\omega)x_0\in{\rm int}N(\theta_{t_0}\omega)$. Noting
that $\hat R(\omega)$ is an invariant random compact set, we have
$\varphi(t_0,\omega)x_0\in\hat R(\theta_{t_0}\omega)$. This is a
contradiction to the fact $\hat R(\omega)\cap{\rm
int}N(\omega)=\emptyset$ for each $\omega$. Therefore we have
obtained $R(\omega)=\hat R(\omega)$, i.e. $R(\omega)$ is a random
repeller with a fundamental neighborhood $\hat N(\omega)$. We now
show $B(R)(\omega)=X-A(\omega)$. In fact we have
\begin{align*}
B(R)(\omega)&=\bigcup_{n\in\mathbb
N}\varphi(n,\theta_{-n}\omega)\hat
N(\theta_{-n}\omega)\\
&=\lim_{n\rightarrow\infty}\varphi(n,\theta_{-n}\omega)\hat
N(\theta_{-n}\omega)\\
&=\lim_{n\rightarrow\infty}\varphi(n,\theta_{-n}\omega)[{\rm int}
N(\theta_{-n}\omega)]^c\\
&=\lim_{n\rightarrow\infty}[\varphi(n,\theta_{-n}\omega){\rm int}
N(\theta_{-n}\omega)]^c\\
&=X-A(\omega),
\end{align*}
where the 4th ``=" follows from the fact that $\varphi(n,\omega)$ is
a homeomorphism on $X$. This completes the proof of the lemma.
\hfill$\Box$

Now we can give the definition of attractor-repeller pair of
$\varphi$.

\begin{definition}\rm
Assume $\varphi$ is an RDS on a compact metric space $X$,
$A(\omega)$ is an attractor of $\varphi$ with a fundamental
neighborhood $N(\omega)$ and the basin of attraction $B(A)(\omega)$.
Then the random set given by
\[
R(\omega)=X\backslash B(A)(\omega)
\]
is called the {\em repeller corresponding to} $A(\omega)$ with the
basin of repulsion $X\backslash A(\omega)$ and a fundamental
neighborhood $X\backslash{\rm int}N(\omega)$. And we call (A,R) an
{\em attractor-repeller pair of} $\varphi$.
\end{definition}

\begin{lemma}\label{attr}
Assume $A(\omega)$ is an attractor with a fundamental neighborhood
$N(\omega)$ and the basin of attraction $B(A)(\omega)$. Then for
arbitrary random closed set $K(\omega)\subset B(A)(\omega)$, there
exists $T(K,\omega)\ge 0$ such that
\begin{equation}\label{ent}
\varphi(t,\omega)K(\omega)\subset {\rm
int}N(\theta_t\omega),~\forall t\ge T(K,\omega).
\end{equation}
\end{lemma}
\noindent{\bf Proof.} For given $K(\omega)\subset B(A)(\omega)$ and
$\forall x\in K(\omega)$, there exists a $t(x)\ge 0$ such that
\[
\varphi(s,\omega)x\in {\rm int}N(\theta_s\omega),~\forall s\ge t(x)
\]
by the definition of basin of attraction and the forward invariance
of ${\rm int}N$ (the forward invariance of ${\rm int}N$ follows from
the fact that $N$ is forward invariant, see page 35 of \cite{Ar1}).
Since $\varphi(s,\omega)$ is a homeomorphism of $X$, there exists an
open neighborhood $U(x)$ of $x$ such that
\begin{equation}\label{enter}
\varphi(s,\omega)U(x)\subset {\rm int}N(\theta_s\omega).
\end{equation}
By the compactness of $K(\omega)$, there exists a finite collection
of such neighborhoods $\{U_i\}_{i=1}^{n}$ which constitutes a finite
open covering of $K(\omega)$ such that (\ref{enter}) hold with $U_i$
instead of $U(x)$. Denote $t_i$ the entrance time of $U_i$ into $N$
and let $T(K,\omega)=\max\{t_i|~i=1,\ldots,n\}$, then we obtain that
(\ref{ent}) holds by the forward invariance of ${\rm int}N$. This
completes the proof of the lemma. \hfill$\Box$

In contrast to the weak attraction in its basin in \cite{Cra}, our
attractor pull-back attracts random closed sets in its basin. See
the following lemma. We remark that the similar result also holds
for repellers and the proof is completely similar.
\begin{lemma}\label{uniatt}
Assume $A(\omega)$ is a random attractor and $B(A)(\omega)$ is the
corresponding basin of attraction, then for any random closed set
$D(\omega)\subset B(A)(\omega)$ we have $A(\omega)$ pull-back
attracts $D(\omega)$.
\end{lemma}
\noindent{\bf Proof.} Assume $N(\omega)$ is a fundamental
neighborhood of $A(\omega)$, then by Lemma \ref{attr} we know that
for any random closed (hence compact, for $X$ being compact) set
$D(\omega)\subset B(A)(\omega)$ there exists a $T_D(\omega)\ge 0$
such that
\[
\varphi(t,\omega)D(\omega)\subset {\rm
int}N(\theta_t\omega),~\forall t\ge T_D(\omega).
\]
For arbitrary non-random $k\in\mathbb N$, by the measure preserving
of $\theta_t$ we obtain that
\begin{align*}
& \mathbb P\{\omega|~
\varphi(t,\theta_{-t}\circ\theta_{-k}\omega)D(\theta_{-t}\circ\theta_{-k}\omega)
\subset N(\theta_{-k}\omega), t\ge T_D(\theta_{-k}\omega)\}\\
=& \mathbb P
\{\omega|~\varphi(t,\theta_{-k}\omega)D(\theta_{-k}\omega)\subset
N(\theta_{t}\circ\theta_{-k}\omega), t\ge T_D(\theta_{-k}\omega)\}\\
=& 1.
\end{align*}
Hence
\begin{align*}
\mathbb P\{\omega|~
&\varphi(k,\theta_{-k}\omega)\circ\varphi(t,\theta_{-t}\circ\theta_{-k}\omega)D(\theta_{-t}\circ\theta_{-k}\omega)\\
&\quad\subset \varphi(k,\theta_{-k}\omega)N(\theta_{-k}\omega), t\ge
T_D(\theta_{-k}\omega)\}=1
\end{align*}
Therefore we have
\[
\overline{\bigcup_{t\ge
T_D(\theta_{-k}\omega)}\varphi(t+k,\theta_{-t-k}\omega)D(\theta_{-t-k}\omega)}
\subset\varphi(k,\theta_{-k}\omega)N(\theta_{-k}\omega)
\]
almost surely. By the definition of omega-limit sets we then obtain
that
\[
\Omega_D(\omega)\subset
\varphi(k,\theta_{-k}\omega)N(\theta_{-k}\omega), \forall
k\in\mathbb N
\]
almost surely and hence
\[
\Omega_D(\omega)\subset \bigcap_{k\in\mathbb
N}\varphi(k,\theta_{-k}\omega)N(\theta_{-k}\omega)=\Omega_N(\omega)=A(\omega)
\]
almost surely. Where $\bigcap_{k\in\mathbb
N}\varphi(k,\theta_{-k}\omega)N(\theta_{-k}\omega)=\Omega_N(\omega)$
holds because $N(\omega)$ is forward invariant, which follows that
\[
\varphi(t,\theta_{-t}\omega)N(\theta_{-t}\omega)\subset\varphi(s,\theta_{-s}\omega)N(\theta_{-s}\omega),
\forall t>s.
\]
Hence we have $A(\omega)$ pull-back attracts $D(\omega)$ by Remark
\ref{re}. This completes the proof of the lemma. \hfill$\Box$

\begin{remark}\rm
With respect to the relation between our definition of attractor and
that of \cite{Cra}, it seems that our definition is stronger. But if
in their Definition 4.1, the fundamental neighborhood is not exactly
the basin of attraction (note that in their definition, basin of
attraction is a special fundamental neighborhood) and the basin
contains the closure of a fundamental neighborhood, then their
definition is equivalent to ours. Besides this, we do not know how
to construct a Lypunov function for their attractor if no further
condition is assumed.
\end{remark}

\begin{remark}\label{sk}\rm
By Lemma \ref{uniatt} we know that an attractor pull-back attracts
any random closed sets inside its basin, but it can not pull-back
attracts its basin itself, for
$\Omega_{B(A)}(\omega)=\overline{B(A)(\omega)}$ by the invariance of
$B(A)(\omega)$. Hence given an attractor $A(\omega)$ (here to
distinguish we call ``attractor" in our definition and call ``weak
attractor" in \cite{Cra}) and an invariant random open neighborhood
$U(\omega)$ of $A(\omega)$ with the property that $A(\omega)$
pull-back attracts any random closed set inside $U(\omega)$, then
$U(\omega)$ must be the basin of attraction of $A(\omega)$. In fact,
if we only know that the attractor $A(\omega)$ attracts any random
closed set inside $U(\omega)$ in probability, the result also holds.
Since in this case, $A(\omega)$ is also a weak attractor defined in
\cite{Cra} and $U(\omega)$ is the basin of it, see Lemma 4.2 and
Proposition 5.1 of \cite{Cra}. But when $A(\omega)$ is regarded as
an attractor, the basin of it should be the same as when it is
regarded as a weak attractor, for the basin being unique.
\end{remark}

\begin{lemma}\label{lya1}
Assume $(A,R)$ is an attractor-repeller pair of $\varphi$, then
there exists an $\mathscr F\times\mathscr B(X)$-measurable Lyapunov function $L$ for $(A,R)$ such that:\\
(i) $L(\omega,x)=0$ when $x\in A(\omega)$, and $L(\omega,x)=1$
     when $x\in R(\omega)$;\\
(ii) for $x\in X\backslash(A(\omega)\bigcup R(\omega))$ and $t>0$,
$1>L(\omega,x)>L(\theta_t\omega,\varphi(t,\omega)x)>0$.
\end{lemma}
\noindent{\bf Proof.} The idea of the proof is originated from
\cite{BS,Ar2}. Assume $N(\omega)$ is a fundamental neighborhood of
$A(\omega)$, and we define the first entrance time of
$\varphi(t,\omega)x$ into $N(\theta_t\omega)$ as follows:
\begin{equation}
\tau(\omega,x):=\left\{
\begin{array}{ll}
    -\infty, & x\in A(\omega); \\
    \inf\{t\in\mathbb R|~\varphi(t,\omega)x\in N(\theta_t\omega)\},
      &  x\in X\backslash(A(\omega)\bigcup R(\omega));\\
    +\infty, & x\in R(\omega). \\
\end{array}
\right.
\end{equation}

Since $\omega\mapsto d(x,N(\omega))$ is measurable, $x\mapsto
d(x,N(\omega))$ is continuous, we have $(\omega,x)\mapsto
d(x,N(\omega))$ is measurable. Hence for arbitrary $t\in\mathbb R$,
$(\omega,x)\mapsto d(\varphi(t,\omega)x,N(\theta_t\omega))$ is
measurable. For $\forall a\in\mathbb R$,
\[
\{(\omega,x)|~\tau(\omega,x)\ge a\}=\bigcap_{t<a,t\in\mathbb
Q}\{(\omega,x)|~d(\varphi(t,\omega)x,N(\theta_t\omega))>0\},
\]
which verifies that $(\omega,x)\mapsto \tau(\omega,x)$ is
measurable.

By the definition of $\tau(\omega,x)$, we have
\begin{align*}
\tau(\theta_t\omega,\varphi(t,\omega)x)&=\inf\{s\in\mathbb
R|~\varphi(s,\theta_t\omega)\circ\varphi(t,\omega)x\in N(\theta_{t+s}\omega)\}\\
&=\inf\{s\in\mathbb
R|~\varphi(t+s,\omega)x\in N(\theta_{t+s}\omega)\}\\
&=\tau(\omega,x)-t.
\end{align*}
Define
\[
L(\omega,x)=\left\{
\begin{array}{ll}
    \frac12 {\rm e}^{\tau(\omega,x)}, & -\infty\le\tau(\omega,x)<0; \\
    \frac12(1+\frac2\pi\arctan\tau(\omega,x)), & 0\le\tau(\omega,x)\le +\infty. \\
\end{array}
\right.
\]
Since $\tau(\omega,x)$ is $\mathscr F\times\mathcal
B(X)$-measurable, hence $L(\omega,x)$ is. It is obvious that the so
defined $L(\omega,x)$ satisfies (i) of Lemma \ref{lya1} and (ii)
follows from the fact
$\tau(\theta_t\omega,\varphi(t,\omega)x)=\tau(\omega,x)-t$. This
terminates the proof of the lemma. \hfill$\Box$

\begin{lemma}\label{strong}
Assume $(A,R)$ is an attractor-repeller pair of $\varphi$ and there
exists a strong fundamental neighborhood $N(\omega)$ of $A(\omega)$,
then there exists a Lyapunov function $L$ for $(A,R)$ with
properties stated in Lemma \ref{lya1} and that $x\mapsto
L(\omega,x)$ is continuous for each $\omega\in\Omega$.
\end{lemma}
\noindent{\bf Proof.} We only need to prove the continuity of
$x\mapsto \tau(\omega,x)$. The proof is completely similar to
Proposition 6.6 of \cite{Ar2}, which is in turn originated from its
deterministic case, see page 71 of \cite{BS}. So we omit details
here. \hfill$\Box$

To distinguish, we call the Lyapunov function obtained in Lemma
\ref{lya1} {\em measurable Lyapunov function} and the one in Lemma
\ref{strong} {\em continuous Lyapunov function}. In contrast to
Lemma \ref{strong} we have the following result.

\begin{lemma}\label{att-rep}
Assume $A(\omega),R(\omega)$ are two disjoint invariant random
compact sets and $L$ is a continuous Lyapunov function for $(A,R)$
with properties stated in Lemma \ref{strong}. Then $(A,R)$ is an
attractor-repeller pair of $\varphi$.
\end{lemma}
\noindent{\bf Proof.} Denote
\[
M(\omega):=\{x|~L(\omega,x)<1\},
\]
then it is easy to see that $R(\omega)=M^c(\omega)$ and hence
$M(\omega)$ is an invariant random open set. For
$\forall0<\alpha<1$, denote
\[
M_\alpha(\omega)=\{x|~L(\omega,x)\le\alpha\}.
\]
Since for any $(x,\omega)\in X\times\Omega$, we have
\[
L(\omega,x)\ge L(\theta_t\omega,\varphi(t,\omega)x),~t\ge 0,
\]
hence $x\in M_\alpha(\omega)$ implies $\varphi(t,\omega,x)\in
M_\alpha(\theta_t\omega)$, i.e.  $M_\alpha(\omega)$ is a forward
invariant random compact set and it is a random neighborhood of
$A(\omega)$. Define
\[
A_\alpha(\omega):=\Omega_{M_\alpha}(\omega)=\bigcap_{T\ge
0}\overline{\bigcup_{t\ge
T}\varphi(t,\theta_{-t}\omega)M_\alpha(\theta_{-t}\omega)},
\]
then by the forward invariance of $M_\alpha$ we have
\[
A_\alpha(\omega)=\bigcap_{t\ge
0}\varphi(t,\theta_{-t}\omega)M_\alpha(\theta_{-t}\omega).
\]
On one hand, we have
\[
A(\omega)=\bigcap_{t\ge
0}\varphi(t,\theta_{-t}\omega)A(\theta_{-t}\omega)\subset\bigcap_{t\ge
0}\varphi(t,\theta_{-t}\omega)M_\alpha(\theta_{-t}\omega)=A_\alpha(\omega).
\]
On the other hand we also have $A_\alpha(\omega)\subset A(\omega)$.
In fact, consider
\[
L(\omega):=\sup_{x\in A_\alpha(\omega)}L(\omega,x).
\]
If the assertion is false, similar to the argument of Proposition
6.2 in \cite{Ar2}, then we have $L(\cdot)>0$ with positive
probability and hence
\[
L(\cdot)>L(\theta_t\cdot),~\forall t>0
\]
with positive probability, a contradiction to the invariance of
$\mathbb P$. Hence we have got that $A=A_\alpha$. Therefore we
obtain that $A(\omega)$ pull-back attracts $M_\alpha(\omega)$ (since
$A_\alpha(\omega)$ does so by its definition), i.e. $A(\omega)$ is
an attractor with $M_\alpha(\omega)$ a fundamental neighborhood. We
now only need to show that $M(\omega)$ is in fact the basin of
attraction of $A(\omega)$, i.e. $B(A)(\omega)=M(\omega)$.

For any random closed set $D(\omega)\subset M(\omega)$ and
$\forall\epsilon>0$, there exists $\alpha<1$ such that
\begin{equation}\label{epsilon}
\mathbb P\{\omega|~D(\omega)\subset M_\alpha(\omega)\}\ge
1-\epsilon.
\end{equation}

By the triangle inequality, we have
\begin{align*}
d(\varphi(t,\omega)D(\omega)|A(\theta_t\omega))\le &
d(\varphi(t,\omega)D(\omega)|\varphi(t,\omega)M_\alpha(\omega))\\
&+d(\varphi(t,\omega)M_\alpha(\omega)|A_\alpha(\theta_t\omega))
+d(A_\alpha(\theta_t\omega)|A(\theta_t\omega)).
\end{align*}
This together with (\ref{epsilon}), the facts $A_\alpha$ attracts
$M_\alpha$ and $A=A_\alpha$ verifies that
\begin{equation}\label{max}
\mathbb
P-\lim_{t\rightarrow\infty}d(\varphi(t,\cdot)D(\cdot)|A(\theta_t\cdot))=0,
\end{equation}
i.e. $A(\omega)$ attracts $D(\omega)$ in probability.

By the above argument and Lemma 4.2 and Proposition 5.1 of
\cite{Cra}, we know that $A(\omega)$ is a weak attractor (defined in
\cite{Cra}) and $M(\omega)$ is the corresponding basin of
attraction. Hence by Remark \ref{sk} we obtain that $M(\omega)$ is
also the basin of $A(\omega)$ when $A(\omega)$ is regarded as an
attractor (defined in present paper). Therefore
$R(\omega)=M^c(\omega)$ is the repeller corresponding to
$A(\omega)$. Hence $(A,R)$ is an attractor-repeller pair of
$\varphi$. This completes the proof of the lemma. \hfill$\Box$

By Lemmas \ref{strong}  and  \ref{att-rep}, we obtain the following
result.

\begin{theorem}\label{th1}
Assume $\varphi$ is an RDS on a compact metric space $X$ and $A,R$
are two disjoint invariant random compact sets. Then $(A,R)$ is an
attractor-repeller pair with strong fundamental neighborhood if and
only if there exists a Lyapunov function $L:\Omega\times
X\rightarrow [0,1]$ such that:\\
(i) $\omega\mapsto L(\omega,x)$ is measurable for each $x\in X$, and
$x\mapsto L(\omega,x)$ is continuous for each
$\omega\in\Omega$;\\
(ii) $L(\omega,x)=0$ when $x\in A(\omega)$, and $L(\omega,x)=1$
     when $x\in R(\omega)$;\\
(iii) for $x\in X\backslash(A(\omega)\bigcup R(\omega))$ and $t>0$,
$1>L(\omega,x)>L(\theta_t\omega,\varphi(t,\omega)x)>0$.
\end{theorem}

\section{Morse decomposition and Lyapunov function}

In this section, we mainly consider the relation between Morse
decomposition and Lyapunov function for RDS.

First, we give the definition of Morse decomposition for random
dynamical systems, which was introduced in \cite{Cra}. For the
deterministic case of Morse decomposition, one can refer to
\cite{Con}.

\begin{definition}\rm (Morse decomposition) Let $\varphi$ be an
RDS on a compact metric space $X$. Assume that $(A_i,R_i)$ are
attractor-repeller pairs of $\varphi$ with
\[
\emptyset=A_0\varsubsetneq A_1\varsubsetneq\cdots\varsubsetneq
A_n=X~{\rm and}~X=R_0\varsupsetneq
R_1\varsupsetneq\cdots\varsupsetneq R_n=\emptyset.
\]
Then the family $D=\{M_i\}_{i=1}^{n}$ of invariant random compact
sets of $X$, defined by
\[
M_i=A_i\bigcap R_{i-1},~1\le i\le n
\]
is called a {\em Morse decomposition} for $\varphi$ on $X$, and each
$M_i$ is called {\em Morse set}. If $D$ is a Morse decomposition,
$M(D)$ is defined to be $\bigcup_{i=1}^nM_i$.
\end{definition}

\begin{remark}\rm
By the definitions of attractor-repeller pair and Morse
decomposition, it is easy to see that $\{\emptyset,X\}$ and
$\{X,\emptyset\}$ are two trivial attractor-repeller pairs and hence
$\{X\}$ is a trivial Morse decomposition for $\varphi$ on $X$.
Moreover, attractor-repeller pair decomposition is a special case of
Morse decomposition. That is, if $(A,R)$ is an attractor-repeller
pair, then $\{M_1(=A),M_2(=R)\}$ is a Morse decomposition.
Conversely, if $D=\{M_i\}_{i=1}^{n}$ is a Morse decomposition, then
we can easily obtain attractor-repeller pairs from it. Moreover,
Morse decompositions can be coarsened. For example, assume
$D=\{M_i\}_{i=1}^{n}$ is a Morse decomposition, where $M_i=A_i\cap
R_{i-1}$. Let $\{i_j\}_{j=1}^{k-1}\subset\{i\}_{i=1}^{n-1}$ and
denote
\[
\tilde A_0=\emptyset,~ \tilde A_j=A_{i_j},~ \tilde A_k=X, ~{\rm
where}~j=1,\cdots,k-1,
\]
then we obtain a coarsened Morse decomposition
\[
\tilde D=\{\tilde M_j\}_{j=1}^{k},~{\rm where}~\tilde M_j=\tilde
A_j\cap\tilde R_{j-1}.
\]
In particular, when $k=2$, the coarsened Morse decomposition
$(\tilde M_1,\tilde M_2)$ is in fact a non-trivial
attractor-repeller pair. It is obvious that we have
\[
M(D)\subset M(\tilde D).
\]
Clearly $D=\{X\}$ is the coarsest Morse decomposition (not
decomposing at all), but there is no finest Morse decomposition, see
Example 2.16 of \cite{Mic} for a deterministic example.
\end{remark}

Similar to the deterministic case, we have the following result
about Morse decomposition for RDS.

\begin{lemma}\label{Mor}
Assume $D=\{M_i\}_{i=1}^{n}$ is a Morse decomposition for $\varphi$
on $X$, then we have
\[
M(D)=\bigcap_{i=0}^n(A_i\bigcup R_i).
\]
\end{lemma}
\noindent{\bf Proof.} The proof is completely similar to that of
Lemma 6 in \cite{Hua}, so we omit the details here. \hfill$\Box$

\begin{lemma}\label{M2}
Assume $D=\{M_1,M_2,\ldots,M_n\}$ is a Morse decomposition for
$\varphi$ on $X$. Then there exists an $\mathscr F\times\mathscr
B(X)$-measurable Lyapunov function $L:\Omega\times
X\rightarrow [0,1]$ such that:\\
(i) $L$ is constant on each $M_i$, i.e. for $\forall x,y\in
M_i(\omega)$, $L(\omega,x)=L(\omega,y)=\alpha_i$, and $\alpha_i$
is independent of $\omega$, $i=1,\ldots,n$;\\
(ii) $\alpha_1<\alpha_1<\cdots<\alpha_n$, i.e.
$L(\cdot,M_1(\cdot))<L(\cdot,M_2(\cdot))<\cdots<L(\cdot,M_n(\cdot))$;\\
(iii) for $x\in X\backslash(\bigcup_{i=1}^{n}M_i(\omega))$ and
$t>0$, $L(\omega,x)>L(\theta_t\omega,\varphi(t,\omega)x)$.
\end{lemma}
\noindent{\bf Proof.} Assume the Morse decomposition
$D=\{M_i\}_{i=1}^{n}$ is determined by attractor-repeller pairs
$(A_i,R_i)$, $i=0,1,\ldots,n$ and assume $l_i(\omega,x)$ is the
Lyapunov function constructed in Lemma \ref{lya1} for the
attractor-repeller pair $(A_i,R_i)$. Let
\begin{equation}\label{ly}
L(\omega,x)=\sum_{i=0}^n\frac{2l_i(\omega,x)}{3^{i+1}},
\end{equation}
then $L(\omega,x)$ is the Lyapunov function desired. In fact, for
the Morse set $M_i(\omega)$, $1\le i\le n$, it is easy to see that
\[
M_i(\omega)\subset A_j(\omega), j\ge i {\rm ~and ~}
M_i(\omega)\subset R_j(\omega), j\le i-1.
\]
Hence by the definition of $l_i(\omega,x)$, we have
$L(\omega,M_i(\omega))=\sum_{j=0}^{i-1}\frac{2}{3^{j+1}}$, which
verifies (i)--(ii) in Lemma \ref{M2}. For $x\in X\backslash
M_D(\omega)$, by Lemma \ref{Mor} we know that there exists an $0\le
i\le n$ such that $x\notin A_i(\omega)\cup R_i(\omega)$. Therefore
we have $l_i(\omega,x)>l_i(\theta_t\omega,\varphi(t,\omega)x)$ for
$\forall t>0$, which together with the fact $l_j(\omega,x)\ge
l_j(\theta_t\omega,\varphi(t,\omega)x)$ for each $0\le j\le n$
verify (iii). \hfill$\Box$

\begin{remark}\label{rem}\rm
If for each Morse set $M_i=A_i\cap R_{i-1}$ in Lemma \ref{M2} $A_i$
has a strong fundamental neighborhood $N_i$, then by Lemma
\ref{strong} and (\ref{ly}) the Lyapunov function obtained in Lemma
\ref{M2} is continuous, i.e. $x\mapsto L(\omega,x)$ is continuous
for each $\omega\in\Omega$.
\end{remark}

\begin{lemma}\label{M1}
Let $D=\{M_1,M_2,\ldots,M_n\}$ be a finite collection of mutually
disjoint invariant random compact sets and assume there exists a
continuous Lyapunov function for $D$ with properties stated in
Remark \ref{rem}, then $D$ is a Morse decomposition for $\varphi$ on
$X$.
\end{lemma}
\noindent{\bf Proof.} Assume $L(\omega,x)$ is a Lyapunov function
for $D$. For definiteness, let $L(\omega,M_i(\omega))=\alpha_i$. By
property (i), (ii) of Lemma \ref{M2}, $\alpha_i$ are non-random
constants and $\alpha_1<\alpha_2<\cdots<\alpha_n$. Let $A_1:=M_1$.
For arbitrary $\alpha_{1,2}$ with $\alpha_1<\alpha_{1,2}<\alpha_2$,
define
\[
N_{1,2}(\omega)=\{x|~\alpha_1\le L(\omega,x)\le\alpha_{1,2}\}.
\]
Then completely similar to the proof of Lemma \ref{att-rep} we know
that $A_1(=M_1)$ is an attractor with a fundamental neighborhood
$N_{1,2}(\omega)$ and the corresponding basin of attraction is
\[
B(A_1)(\omega)=\{x|~\alpha_1\le L(\omega,x)<\alpha_{2}\}.
\]
Therefore the repeller $R_1$ corresponding to $A_1$ is
\[
R_1(\omega)=\{x|~L(\omega,x)\ge\alpha_{2}\}.
\]
Hence $M_2,\cdots$, $M_n\subset R_1$.

For $\forall\alpha_{2,3}\in(\alpha_2,\alpha_3)$, define
\[
N_{2,3}(\omega)=\{x|~\alpha_1\le L(\omega,x)\le\alpha_{2,3}\}.
\]
It is obvious that $M_1\bigcup M_2\subset N_{2,3}$ and $N_{2,3}$ is
a fundamental neighborhood. Assume $A_2$ is the attractor inside
$N_{2,3}$, i.e.
\begin{equation}\label{A2}
A_2(\omega)=\bigcap_{t\ge0}\varphi(t,\theta_{-t}\omega)N_{2,3}(\theta_{-t}\omega).
\end{equation}
Hence we have $M_1\bigcup M_2\subset A_2$. Therefore we have
obtained $A_2\bigcap R_1\supset M_2$, next we show that $A_2\bigcap
R_1\subset M_2$. Since for any $x\in
N_{2,3}(\omega)\backslash(M_1(\omega)\bigcup M_2(\omega))$ and
$\forall t>0$, we have
\[
L(\theta_t\omega,\varphi(t,\omega)x)<L(\omega,x).
\]
Therefore, by the proof of Lemma \ref{att-rep}, for
$\forall\alpha\in(\alpha_2,\alpha_3)$, the forward invariant random
compact set
\[
N_\alpha(\omega)=\{x|~\alpha_1\le L(\omega,x)\le\alpha\}
\]
is always a fundamental neighborhood of $A_2(\omega)$. Hence we have
\[
A_2(\omega)\subset\bigcap_{n\in\mathbb
N}N_{\alpha_2+\frac1n}(\omega),
\]
and similarly we also have
\[
R_1(\omega)\subset\bigcap_{n\in\mathbb N}\tilde
N_{\alpha_2-\frac1n}(\omega),
\]
where
\[
N_{\alpha_2+\frac1n}(\omega)=\{x|~\alpha_1\le
L(\omega,x)\le\alpha_2+\frac1n\},~\tilde
N_{\alpha_2-\frac1n}(\omega)=\{x|~L(\omega,x)\ge\alpha_2-\frac1n\}.
\]
Thus
\begin{align*}
A_2(\omega)\bigcap R_1(\omega)&\subset(\bigcap_{n\in\mathbb
N}N_{\alpha_2+\frac1n}(\omega))\bigcap(\bigcap_{n\in\mathbb
N}\tilde N_{\alpha_2-\frac1n}(\omega)) \\
&\subset\bigcap_{n\in\mathbb N}
(N_{\alpha_2+\frac1n}(\omega)\bigcap\tilde
N_{\alpha_2-\frac1n}(\omega)) \\
&=\{x|~L(\omega,x)=\alpha_2\}=M_2(\omega),
\end{align*}
i.e. we have obtained $A_2\bigcap R_1=M_2$. Then we can obtain $R_2$
from $A_2$, i.e.
\[
R_2(\omega)=\{x|~L(\omega,x)\ge\alpha_{3}\}.
\]

Similar to the above arguments, let
\[
N_{3,4}(\omega)=\{x|~\alpha_{1}\le
L(\omega,x)\le\alpha_{3,4}\},~{\rm
where}~\alpha_{3,4}\in(\alpha_3,\alpha_4),
\]
and we immediately obtain $A_3$ similar to (\ref{A2}). Hence we at
once obtain the repeller $R_3$ corresponding to $A_3$. Inductively,
we can obtain $A_4$, $R_4$, $\ldots$, $A_{n-1}$, $R_{n-1}$ in the
same way. Let $A_0=R_n=\emptyset,A_n=R_0=X$. Therefore we have
obtained
\[
\emptyset=A_0\varsubsetneq A_1\varsubsetneq\cdots\varsubsetneq
A_n=X~{\rm and}~X=R_0\varsupsetneq
R_1\varsupsetneq\cdots\varsupsetneq R_n=\emptyset
\]
from $M_i,i=1,\ldots,n$ satisfying
\[
M_i=A_i\bigcap R_{i-1},~1\le i\le n.
\]
This shows that $D$ is a Morse decomposition for $\varphi$ on $X$
and hence completes the proof of the lemma. \hfill$\Box$

By Lemmas \ref{M2}, \ref{M1} and Remark \ref{rem} we obtain the
following theorem.

\begin{theorem}\label{th2}
Assume $\varphi$ is an RDS on a compact metric space $X$ and let
$D=\{M_1,M_2,\ldots,M_n\}$ be a finite collection of mutually
disjoint invariant random compact sets. Then $D$ is a Morse
decomposition for $\varphi$ on $X$ with each $A_i$ having a strong
fundamental neighborhood if and only if there exists a Lyapunov
function $L:\Omega\times
X\rightarrow [0,1]$ such that:\\
(i) $\omega\mapsto L(\omega,x)$ is measurable for each $x\in X$, and
$x\mapsto L(\omega,x)$ is continuous for each
$\omega\in\Omega$;\\
(ii) $L$ is constant on each $M_i$, i.e. for $\forall x,y\in
M_i(\omega)$, $L(\omega,x)=L(\omega,y)=\alpha_i$, and $\alpha_i$
is independent of $\omega$, $i=1,\ldots,n$;\\
(iii) $\alpha_1<\alpha_1<\cdots<\alpha_n$, i.e.
$L(\cdot,M_1(\cdot))<L(\cdot,M_2(\cdot))<\cdots<L(\cdot,M_n(\cdot))$;\\
(iv) for $x\in X\backslash(\bigcup_{i=1}^{n}M_i(\omega))$ and $t>0$,
$L(\omega,x)>L(\theta_t\omega,\varphi(t,\omega)x)$.
\end{theorem}

\begin{remark}\rm
By Lemmas \ref{lya1} and \ref{M2}, we can construct measurable
Lyapunov functions for attractor-repeller pairs and Morse
decompositions. But to construct continuous Lyapunov functions, by
Theorems \ref{th1} and \ref{th2} we see that we must find a strong
fundamental neighborhood for a given attractor, which is not an easy
thing. Note that the construction in \cite{Ar2}, which follows from
\cite{BS}, is not applicable when $A$ is not globally attracting.
For deterministic case, an invariant compact set $A$ is called an
{\em attractor} if there exists a fundamental neighborhood $U$ of
$A$ such that the omega-limit set of $U$, $\Omega_U=A$, see
\cite{Con}. This implies that there exists a strong fundamental
neighborhood $U$ of $A$ such that $\Omega_U=A$, see Proposition 1.9
on page 409 of \cite{Rob} for details. But for random case, we do
not know whether or not similar result holds. That is, we do not
know generally how to construct a strong fundamental neighborhood.
Therefore we do not request that the fundamental neighborhood of an
attractor be a strong  one in Definition \ref{def}. Note also that
the construction of Proposition 1.10 on page 409 of \cite{Rob} does
not hold in the random setting. The main difficulty for these stems
form the non-uniformity and the non-autonomy of RDS. This
non-uniformity is one of the essential features of RDS.
\end{remark}

Now we give a simple example to illustrate our results. The example
is borrowed from \cite{Cra}, which is also used in \cite{Liu}.

\begin{example}\rm
Consider the Stratonovich stochastic differential equation (SDE)
\begin{equation}\label{exam}
{\rm d}X_t=(X_t-X_t^3){\rm d}t+(X_t-X_t^3)\circ{\rm d}W_t
\end{equation}
on the interval $[-1,1]$. To put a stochastic differential equation
in the framework of RDS, we model white noise as a metric dynamical
system as follows: Let $\Omega$ be the space of continuous functions
$\omega:\mathbb R\rightarrow\mathbb R$ satisfying that
$\omega(0)=0$, let $\mathscr F$ be the Borel sigma-algebra induced
by the compact-open topology of $\Omega$, and let $\mathbb P$ be the
Wiener measure on $(\Omega,\mathscr F)$, i.e. the distribution on
$\mathscr F$ of a standard Wiener process with two-sided time. The
shift $\theta_t$ is defined by
$\theta_t\omega(s)=\omega(t+s)-\omega(t)$. Then $(\Omega,\mathscr
F,\mathbb P,(\theta_t)_{t\in\mathbb R})$ is an ergodic metric
dynamical system driving the SDE (\ref{exam}), and
$W_t(\omega)=\omega(t)$. See Appendix A.3 of \cite{Ar1} for details.

From p.123 of \cite{Klo} we know that the RDS $\varphi:\mathbb
R\times\Omega\times [-1,1]\mapsto [-1,1]$ generated by SDE
(\ref{exam}) can be expressed by
\[
\varphi(t,\omega)x=\frac{x{\rm e}^{t+W_t(\omega)}}{(1-x^2+x^2{\rm
e}^{2t+2W_t(\omega)})^{\frac12}}.
\]
Hence
\[
\varphi(t,\theta_{-t}\omega)x=\frac{x{\rm
e}^{t-W_{-t}(\omega)}}{(1-x^2+x^2{\rm
e}^{2t-2W_{-t}(\omega)})^{\frac12}}.
\]
Consider the interval $N:=[1/2,1]$, then we can easily see that
\[
\Omega_N(\omega)\equiv\{1\}
\]
by the fact that $\lim_{t\rightarrow\infty}\frac{W_t}{t}=0$ almost
surely. Hence $\{1\}$ is a an attractor with a fundamental
neighborhood $[1/2,1]$. And we can obtain a forward invariant
fundamental neighborhood of $\{1\}$ as
\[
\tilde
N(\omega)=\overline{\bigcup_{t\ge0}\varphi(t,\theta_{-t}\omega)[1/2,1]}.
\]
Clearly the basin of attraction of $\{1\}$ is  $(0,1]$ and hence the
corresponding repeller to it is $[-1,0]$. By Lemma \ref{lya1} there
exists a Lyapunov function $L(\omega,x)$ for the attractor-repeller
pair which is $0$ when $x=1$, is strictly decreasing when $x\in
(0,1]$, is $1$ when $x\in[-1,0]$. Similarly $\{-1\}$ is an attractor
with basin of attraction $[-1,0)$ and the corresponding repeller is
$[0,1]$. Therefore $\{-1,1\}$ is also an attractor with basin of
attraction $[-1,0)\cup (0,1]$ and the corresponding repeller is
$\{0\}$. If we set $A_0=\emptyset, A_1=\{-1\},A_2=\{-1,1\},A_3=X$,
then the corresponding repellers are
$R_0=X,R_1=[0,1],R_2=\{0\},R_3=\emptyset$. Consequently the
corresponding Morse sets are $M_1=\{-1\},M_2=\{1\},M_3=\{0\}$. Thus
by Lemma \ref{M2}, there exists a Lyapunov function for this Morse
decomposition. If we initially set $A_0=\emptyset, A_1=\{1\},
A_2=\{-1,1\}, A_3=X$, then the corresponding Morse sets are
$M_1=\{1\},M_2=\{-1\},M_3=\{0\}$ and we can obtain the similar
result.
\end{example}

\section*{Acknowledgements} The authors express their sincere thanks to Professor Yong Li for
his instructions and many invaluable suggestions.


\begin{thebibliography}{99}

\bibitem{Ar1}
L. Arnold, Random Dynamical Systems, Springer-Verlag, Berlin
Heidelberg New York, 1998.

\bibitem{Ar2}
L. Arnold and B. Schmalfuss, Lyapunov's second method for random
dynamical systems, J. Differential Equations, 177 (2001) 235--265.

\bibitem{BS}
N.P. Bhatia and G.P. Szeg\"o, Stability Theory of Dynamical Systems,
Springer-Verlag, Berlin, 1970.


\bibitem{Cas}
C. Castaing, and M. Valadier,  Convex Analysis and Measurable
Multifunctions. (Lecture Notes in Mathematics, vol 580).
Springer-Verlag, Berlin Heidelberg New York, 1977.

\bibitem{Chu}
I. Chueshov, Monotone Random systems Theory and Applications,
(Lecture Notes in Mathematics vol 1779), Springer-Verlag, Berlin,
2002.

\bibitem{Con}
C. Conley, Isolated Invariant Sets and the Morse Index, Conf. Board
Math. Sci.,  vol 38, Amer. Math. Soc., Providence, 1978.

\bibitem{Cr}
H. Crauel, Random Probability Measures on Polish Spaces (Stochastics
Monographs vol 11), Taylor \& Francis, London, 2002.

\bibitem{Cra}
H. Crauel, L.H. Duc and S. Siegmund, Towrds a Morse theory for
random dynamical systems, Stochastics and Dynamics, 4 (2004)
277-296.

\bibitem{Hua}
T. Huang, Decompositions and Liapounov functions, Chaos, Solitons
and Fractals, 13 (2002) 209-214.

\bibitem{Klo}
P.E. Kloeden and E. Platen, Numerical Solution of Stochastic
Differential Equations, Springer-Verlag, Berlin Heidelberg New York,
1992.

\bibitem{Lia}
M.A. Liapounoff, Probl\`eme g\'en\'eral de la stabilit\'e du
mouvement, Ann. Fac. Sci. Toulouse 9 (1907). [Translation of the
Russian edition, Kharkov 1892, reprinted by Princeton University
Press, Princeton, NJ, 1949 and 1952.]

\bibitem{Liu}
Z. Liu, The random case of Conley's theorem, Nonlinearity 19 (2006)
277-291.

\bibitem{Mic}
K. Mischaikow and M. Mrozek, Conley index, in Handbook of dynamical
systems, Vol. 2, 393--460, North-Holland, Amsterdam, 2002.

\bibitem{Rob}
C. Robinson, Dynamical Systems: Stability, Symbolic Dynamics, and
Chaos, 2nd edition, CRC Press, 1998.


\bibitem{UV}
S.M. Ulam and J. von Neumann, Random ergodic theorems, Bull. Amer.
Math. Soc. 51 (1945) 660.

\end{thebibliography}
\end{document}